\documentclass[11pt,english]{article}
\usepackage{epsfig,graphicx,amsfonts,amsmath,amsthm,amscd,amssymb,latexsym}

\numberwithin{equation}{section}

\newtheorem{theorem}{Theorem}
\newtheorem{lemma}{Lemma}
\newtheorem{corollary}{Corollary}

\begin{document}

\title{\bf
Global in Time Asymptotic Solutions to Kolmogorov--Feller-Type
Parabolic Pseudodifferential Equations with a Small Parameter.
Forward and Backward in Time Motion.
\thanks{This work was supported by DFG project 436 RUS 113/895/0-1.}}
\author{V.~G.~Danilov}

\date{}

\maketitle

\begin{abstract}
The goal of the present paper is to present a new approach
to the construction of asymptotic (approximating) solutions to parabolic PDE
by using the characteristics.
\end{abstract}

\section{\bf Forward in time motion}

The goal of the present paper is to present a new approach
to the construction of asymptotic (approximating) solutions to parabolic PDE
by using the characteristics.
This approach allows one to construct global in time solutions not only
for the usual Cauchy problems but also for the inverse problems.
We will work with Kolmogorov--Feller-type equations with diffusion,
potential, and jump terms. The equation under study has the form:
$$
-\varepsilon\frac{\partial u}{\partial t}
+P\bigg(\overset{2}{x},-\varepsilon\overset{1}{\frac{\partial}{\partial
x}}\bigg)u=0,
$$
$$
u(x,t,\varepsilon)|_{t=0}=e^{-S_0(x)/\varepsilon}\varphi^0(x),
$$
where $P(x,\xi)$ is the symbol of the Kolmogorov--Feller operator,
$\varepsilon\to+0$ is a small parameter characterizing the
frequency and the amplitude of jumps of the Markov stochastic
process with transition probability given by $P(x,\xi)$.
To be more precise, we bear in the mind the following form of $P(x,\xi)$:
\begin{equation}\label{t.0}
P(x,\xi)=(A(x)\xi,\xi) +V(x)+(B(x),\xi)+
\int_{\mathbb{R}^n}\bigg(e^{i (\xi,\nu)}-1\bigg)\mu(x,d\nu),
\end{equation}
where $A(x)$ is a positive smooth matrix, $\mu(x,d\nu)$ is a
family of positive bounded measures smooth with respect to $x$
such that
$$
\int_{\mathbb{R}^n}\nu_i\mu_{x_i}(x,d\nu)=0,\qquad i=1,\dots,n,
$$
and $B$ and $V$ are smooth in $x$ (more precise conditions see
below). The construction of forward in time global asymptotic
solution to equations of this type was developed by V.~Maslov,
\cite{M1,8,10}, for a version of this construction, see also in
\cite{DD,9,D5}. Maslov's approach is based on ideas similar to
those used in his famous canonical operator construction (or in
Fourier integral operators theory). This construction is based
on some integral representation and is not suitable for
constructing backward in time solutions. Another approach to
the global asymptotic solution construction was suggested in
\cite{DAl} and is based on the construction of generalized
solutions to continuity equation in a discontinuous velocity
field~\cite{1}.

We assume that the class of solutions under study admits the following limits:

(1) -- logarithmic pointwise limit $\lim_{\varepsilon \to0}(-\varepsilon \rm{ln}u)$.
We denote this limit by $S=S(x,t)$ and assume that it is a piecewise smooth function
with bounded first-order derivatives and a singular support in the form of a
stratified manifold~$M$.

(2) -- weak limit of the expression $\exp{(2S/\varepsilon)}u^2$.
We denoted it by $\rho$ and assume that $\rho$ is the sum of the
function ($\rho_{\text{reg}}$) smooth outside~$M$ and the Dirac $\delta$-function
on $M$. Note that here we deal with the limit in the weighted weak sense!

If $S$ and $\rho$ are smooth function, then the following representation is true
(in the usual sense)
$$
u=\exp{(-S/\varepsilon)}\sqrt \rho_{\text{reg}}(1+o(\varepsilon)).
$$

{\bf Example 1:}

$$
u|_{t=0}=\exp{(-S_0/\varepsilon)}\varphi^0(x),
$$
where $S_0\geq0$ is a smooth function, $\varphi^0\in C^\infty_0$.

Here a WKB-like approach can be used (Yu.~Kifer,~\cite{K1}; \, V.~Maslov,~\cite{M1,8,10}).
It gives an asymptotic (approximating) solution in the form (cf.~\cite{MF})
$$
u_{\text{as}}=u_{\text{as}}(x,t,\varepsilon)=\exp{(-S(x,t)/\varepsilon)}(\varphi_0(x,t)+\dots
+\varepsilon^k\varphi_k(x,t))
$$
for arbitrary $k$. Here $S(x,t)$ is the solution to the Cauchy problem for the
Hamilton--Jacobi equation
\begin{equation}\label{t.2}
S_t+P(x,\nabla S)=0,
S_{t=0}=S_0(x),
\end{equation}
and $\varphi_0=\varphi_0(x,t)$ is the solution to the transport equation
\begin{equation} \label{t.3}
\varphi_{0t}+(\nabla_{\xi}P(x,\nabla S),\nabla\varphi_0)+
\frac{\varphi_0}{2}{\rm tr}(P_{\xi \xi}(x,\nabla S)S_{xx}=0,
\end{equation}
$$
S_{t=0}=S_0(x).
$$
Both of the solutions $S$ and $\varphi_0$ are defined via solutions of the Hamilton system
\begin{equation}\label{t.4}
\dot x=\nabla_\xi P(x,p),\qquad x|_{t=0}=\alpha,
\dot p=-\nabla_x P(x,p),\qquad p|_{t=0}=\nabla S_0(\alpha).
\end{equation}

They are smooth while
$$
\frac{Dx}{D\alpha}\ne 0.
$$

There are symplectic geometry objects corresponding to this construction:

(1) the phase space $\mathbb{R}^{2n}_{x,p}=\mathbb{R}^n_x \times  \mathbb{R}^n_p$;

(2) the Lagrangian manifold $\Lambda_n^t\in\mathbb{R}^{2n}_{x,p}$,
$$
\Lambda_n^0=(x=\alpha, p=\nabla S_0(\alpha)),
$$
$$
\Lambda_n^t=g^t_P\Lambda_n^0,
$$
where $g^t_P$ is a shift mapping along the Hamiltonian system trajectories;

(3) the projection mapping $\pi: \Lambda_n^t\to \mathbb{R}^n_x$
with Jacobi matrix $\frac{\partial x}{\partial\alpha}$.
\vspace{6mm}

The main assumption that is required is the following one.

The trajectories of the Hamilton system form a manifold of the phase space
(at least in the area of the phase space under study).

Let $Dx/D\alpha\ne 0$ for $t\in [0,T]$, then we have the following statement
(V.~Maslov,~\cite{M1,8}; V.~Danilov, \cite{DD,9,D5}):

\begin{theorem}
The inequality
$$
{\rm Re} P(x,p+i\eta)\le P(x,p),\qquad \eta \in \mathbb{R}^n
$$
is necessary and sufficient for the following estimation to hold:
$$
\|\exp{(S(x,t)/\varepsilon)}(u_{\text{as}}-u)\|_{C(\mathbb{R}^{n_x})}\le C_M\varepsilon^M,
$$
where $M=M(k)\to\infty$ as $k\to \infty$ and the function $S$ is a solution to the
Hamilton--Jacobi equation with Hamiltonian $P(x,p)$.
\end{theorem}

It is easy to verify that the function $P=P(x,\xi)$ -- the symbol introduced above --
satisfies the inequality mentioned in the theorem.

{\bf Example 2:}

$$
-\varepsilon u_t+\varepsilon^2u_{xx}=0, \qquad u|_{t=0}=\exp{(-S_0/\varepsilon)}\varphi^0.
$$

The corresponding system is:
$$
S_t+(S_x)^2=0,  (Hamilton--Jacobi \ equation),
$$
\vspace{1mm}
$$
\varphi_{0t}+2S_x\varphi_{0x}+S_{xx}\varphi_0=0, (transport \  equation)
$$
\vspace{1mm}
$$
\dot x=2p,\qquad x|_{t=0}=\alpha,
$$
$$
\dot p=0,\qquad p|_{t=0}=\frac{\partial S_0}{\partial \alpha}, \qquad (Hamilton \ system).
$$
\vspace{1mm}

Its solution has the form
$$
x=\alpha+2t\frac{\partial S_0}{\partial\alpha},\qquad p=p|_{t=0}
$$
and
$$
Dx/D\alpha=1+2t \frac{\partial^2 S_0}{\partial\alpha^2}.
$$
(i)
\begin{figure}[h!] 
\includegraphics[width=0.4\textwidth]{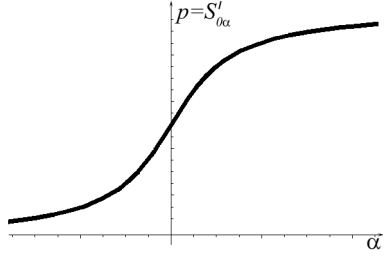}
\caption{ Solution $\alpha=\alpha(x,t)$ exists for all $t$}
\label{fig01}
\end{figure}

(ii)
\begin{figure}[ht!] 
\includegraphics[width=0.4\textwidth]{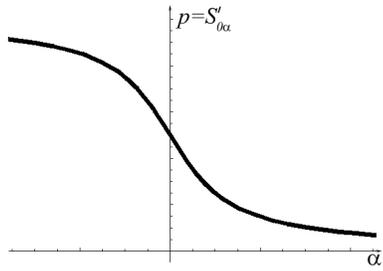}
\caption{Solution $\alpha=\alpha(x,t)$ exists for $t<t^*=\max|2S''_{0\alpha \alpha}|^{-1}$}
\label{fig02}
\end{figure}

For $t>t^*$ in case (ii), we get $\Lambda^1_t$ of the shape plotted in Fig.~3.

\begin{figure}[ht!] 
\includegraphics[width=0.4\textwidth]{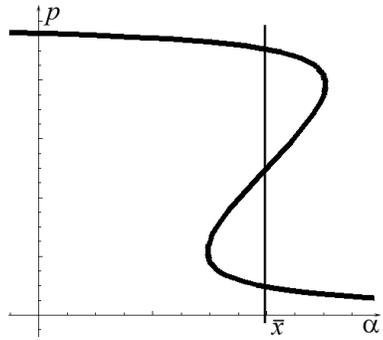}
\caption{Lagrangian curve for $t>t^*$.}
\label{fig03}
\end{figure}

One can see that in the case under study there are three values of $S$ at the
point $\bar{x}$.
This means that we can present an asymptotic solution near this point in the form
of a linear combination:
\begin{equation}\label{d.9}
u=\sum^{3}_{j=1}c_ju_j,
\end{equation}
where each of the functions $u_j=\exp(-S_j/\varepsilon)\varphi_j$, $j=1,2,3$,  satisfies the equation with the same accuracy.
But the functions themselves are not equivalent in contrast to the hyperbolic case.

For example, it is clear that if the inequality
$$
S_1(\bar{x},t)>S_2(\bar{x},t)
$$
holds at a certain point $\bar{x}$, then the ``WKB'' solutions
$u_1$ and $u_2$ at the point $\bar{x}$ satisfy the relation
\begin{align}\label{d.10}
u_1|_{x=\bar{x}}
&=e^{-S_1(\bar{x},t)/\varepsilon}\varphi_1(\bar{x},t)
\nonumber\\
&=e^{-S_1(\bar{x},t)/\varepsilon}\varphi_2(\bar{x},t)
\big(e^{-(S_1-S_2)/\varepsilon}\varphi_1/\varphi_2)|_{x=\bar{x}}
\nonumber\\
&=u_2|_{x=\bar x} O(\varepsilon^N),
\end{align}
where $N>0$ is an arbitrary number.
This follows from the fact that the difference
$(S_1-S_2)|_{\bar{x}}$ in parentheses in the exponent is positive.

Thus, at each point in formula \eqref{d.9}, it is necessary to choose
the term where the function $S_j$ is minimal.
Such a choice leads to an expression of the form
\begin{equation}\label{d.11}
u=e^{-\Phi(x,t)/\varepsilon}\varphi(x,t),
\end{equation}
where $\Phi=\Phi(x,t)=\min_x\{S_j(x,t)\}$.
It is clear that expression \eqref{d.11} is the leading term
of the approximate solution.

The corresponding Lagrangian manifold is the following one:

\begin{figure}[h!] 
\includegraphics[width=7cm]{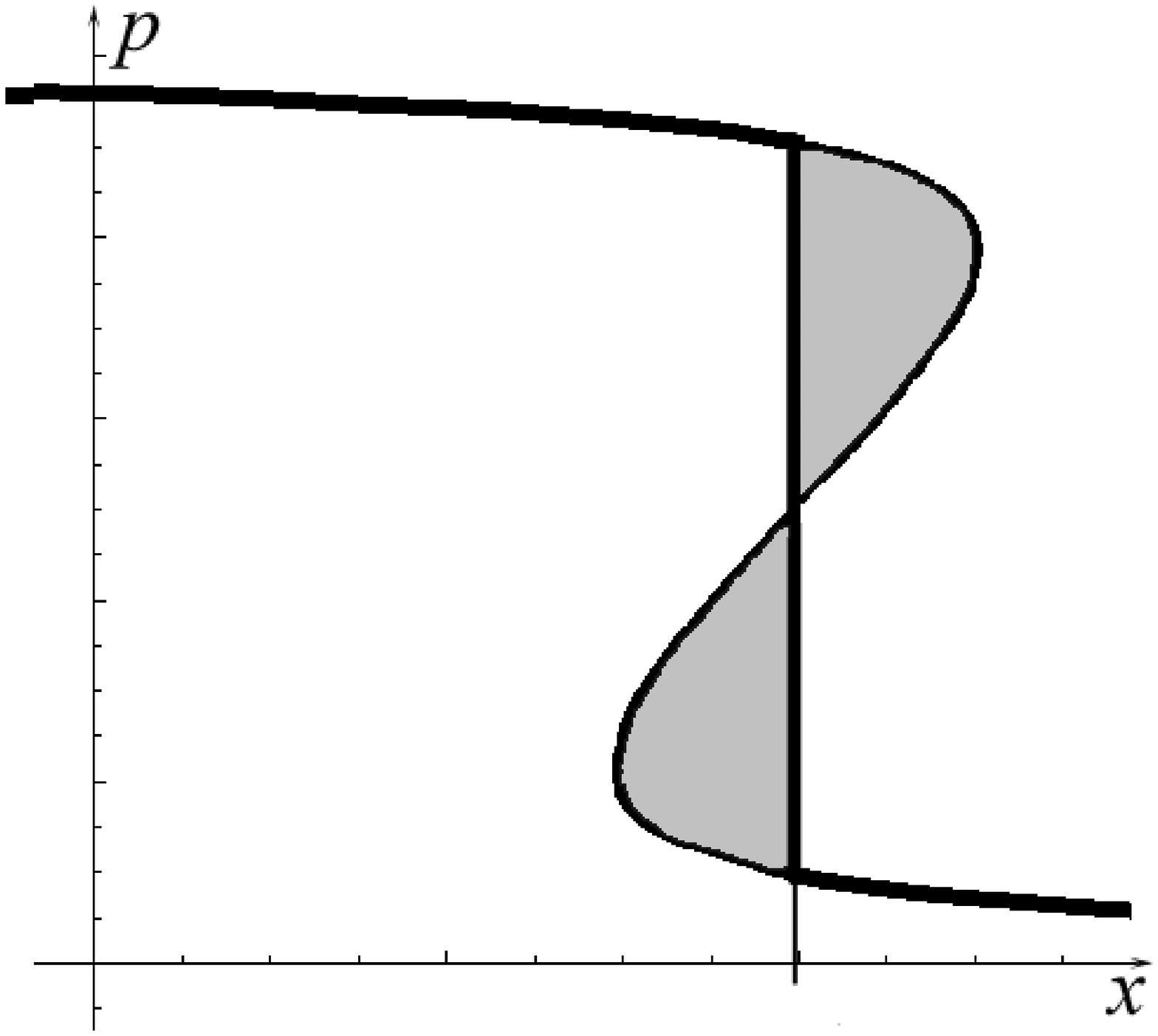}
\caption{}
\end{figure}

The vertical line position is such that the above-mentioned squares are equal.

It is interesting to note that, in the one-dimensional case,
there is a direct connection between the Hamilton--Jacobi equation
$$
S_t+H(x,S_x)=0
$$
and the conservation law of the form
$$
u_t+\frac{\partial}{\partial x}H(x,u)=0,
$$
where $u=S_x$ and the velocity of the vertical line is defined by Rankine--Hugoniot condition
corresponding to the conservation law.
If the Lagrangian manifold has a jump, then the corresponding value function
$$
\Phi(x,t)={\rm {min}}_{x}S(x,t)
$$
has a jump in the first derivative.

Assumption of the actual analyticity for all objects provides that there is no
"concentration of singularities", --- each one can be considered separately
in a sense and one can construct a value function as a solution of the
Hamilton--Jacobi--Bellman equation.

What about the amplitude function? As was mentioned above, it is a solution of
the transport equation
\begin{equation}
\varphi_{0t}+(\dot x,\nabla\varphi_0)+
\frac{\varphi_0}{2}{\rm tr}(P_{\xi \xi}(x,\nabla S)S_{xx}=0,
\end{equation}

Generally, the velocity field calculated from the Hamilton system has a jump
simultaneously with a jump in $p=\nabla S$ (and then in $\dot x$).
Thus, the problem (still open!) is to solve the transport equation in a discontinuous
velocity field. We avoid this problem considering the squared solution
of the transport equation $\rho=\varphi_0^2$.
Madelung,~\cite{Mad} (about 100 years ago!) observed that it satisfies the continuity equation
$$
\rho_t+(\nabla, \dot x \rho)+\frac{\rho}{2}\rm{tr}P_{x\xi}=0.
$$
in the smooth case.

Our case is more complicated: we again have a discontinuity velocity field.
There are few approaches to this problem solution:
the theory of measure solutions (F.~Murat, P.~LeFloh, B.~Hayes, T.~Zang, Y.~Zheng et al.,
\cite{2,3,5,LF,TZ,Y}),
the box approximation (M.~Oberguggenberger, M.~Nedel'kov, \cite{NOb}),
the weak asymptotics method (V.~Danilov, V.~Shelkovich, \cite{15,16,17,18,19}),
the generalized characteristics method (V.~Danilov, D.Mitrovich and V.~Danilov,
\cite{11,12,14,Ddd}).
The last allows one to construct a solution to the continuity equation in the case where
a singular support of the velocity field is a stratified manifold with smooth strata
which are transversal to the (incoming!) trajectories of the velocity field.

If for some time interval $[t_1, t_2]$, the singular support of the velocity field
preserves its structure (the mapping of the singular support induced by the shift
along the Hamilton flow is a diffeomorphism), then it is possible to show
that the singular support of the velocity has the required structure.
If the structure is changing (e.g., a jump appears, see Fig.5 and Fig.4-the
last step of evolution in time),
then one can use the weak asymptotics method to construct a global solution
to the Hamilton--Jacobi and continuity equations. This approach is based
on a  "new (generalized) characteristics" constructed
by V.~Danilov and D.~Mitrovic, \cite{11,12,Ddd} in the case where the strata
of the singular support are of codimension~1.

\begin{figure}[h!] 
\includegraphics[width=\textwidth]{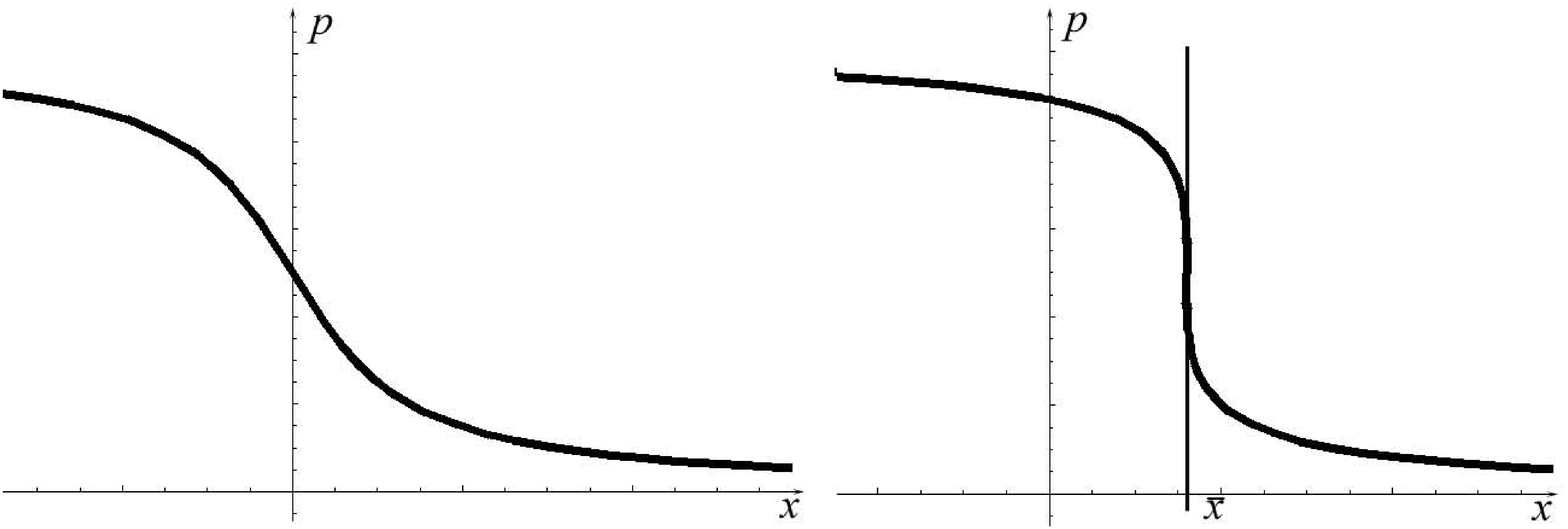}
\caption{}
\end{figure}

\begin{figure}[h!] 
\includegraphics[width=7cm]{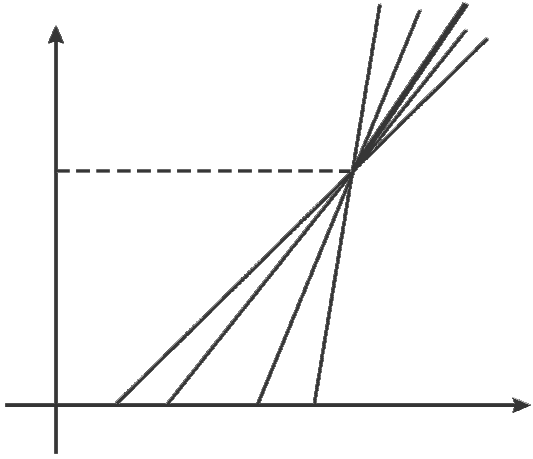}
\caption{Classical characteristics}
\end{figure}

The main idea of this approach is to consider the singularity origination
as a result of nonlinear solitary wave interaction.

\begin{figure}[h!] 
\includegraphics[width=7cm]{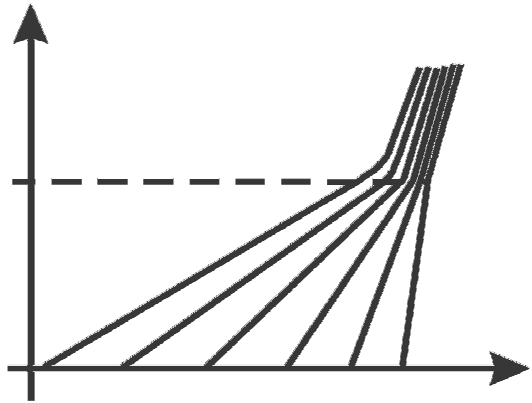}
\caption{New characteristics}
\end{figure}

A simple example is the Hamilton flow corresponding to the heat equation
from the previous example. The Hamilton--Jacobi equation in this case
is equivalent to the Hopf equation for the momentum $p$:
$$
p_t+(p_x)^2=0.
$$

The solution in this case has the form
$$
p=p_0+a(H(\phi_1-x)(\phi_1-x)-H(\phi_2-x)(\phi_2-x))
$$
and is plotted below.

\begin{figure}[h!] 
\includegraphics[width=10cm]{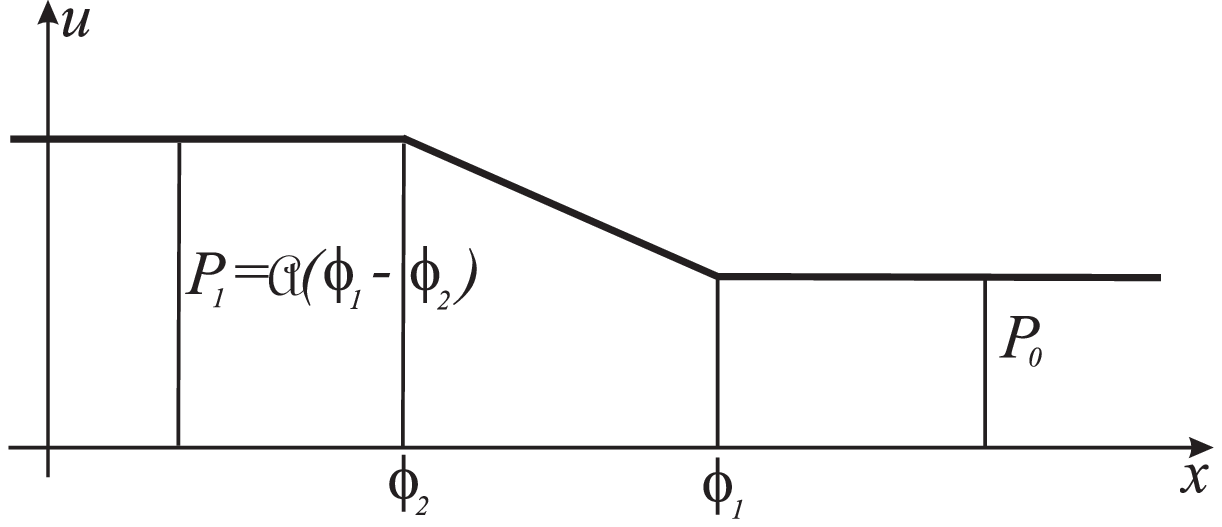}
\caption{Plot of the function $p$}
\end{figure}

To consider $p^2$, we have to calculate the product
$$
H(\phi_1-x)H(\phi_2-x).
$$
Here the following equality holds:
$$
H(\phi_1-x)H(\phi_2-x)=
B(\frac{\phi_2-\phi_1}{\mu})H(\phi_1-x)
$$
$$
+(1-B(\frac{\phi_2-\phi_1}{\mu}))H(\phi_2-x)+O_{D'}(\mu),
$$
where $\mu$ is an arbitrary small parameter and $O_{D'}(\mu)$
is a small quantity in the sense of distributions,
$$
(O_{D'}(\mu),\psi)=O(\mu),
$$
for each $\psi$, which is a test function.
The time evolution of the function $p$ is such that
the slanting intercept of a straight line preserves its shape
until it takes the vertical position and then a jump begins to propagate.
This means that, at every time instant,
the solution anzatz can be presented in the form of a linear combination of
Heaviside functions. This allows one to use a formula which
express the product of Heaviside functions as their linear combination,
and hence, uniformly in time, we see that the functions $p$ and
$p^2$ (the last up to a small quantity) belong to the same linear space,
for detail, see \cite{19,11,12,Ddd}.

Thus, we can prove the following theorem.

\begin{theorem}
Assume  that the following conditions are satisfied
for $t\in[0,T]$, $T>0$:

{\rm(1)} there exists a smooth solution of the Hamiltonian system,

{\rm(2)} the singularities of the velocity field
$$
u=\nabla_\xi P(x,\nabla S)
$$
form a stratified manifold with smooth strata and
$\operatorname{Hess_\xi} P(x,\xi)>0$.

Then there exists a generalized solution $\rho$ of the Cauchy
problem for continuity equation in the sense of the integral
identity introduced in \cite{17,18} and at the points  where
the projection $\pi$ is bijective, the asymptotic solution of
the Cauchy problem for Kolmogorov-Feller type equation has the
form
$$
u=\exp(-S(x,t)/\varepsilon)(\sqrt{\rho_{\text{reg}}}+O(\varepsilon)).
$$
\end{theorem}

\setcounter{equation}{0}

\section{\bf Backward in time motion}

As was shown above, all that we need to go forward in time is the Hamilton system:
\begin{equation}
\dot x=\nabla_\xi P(x,p),\qquad x|_{t=0}=\alpha,
\end{equation}
$$
\dot p=-\nabla_x P(x,p),\qquad p|_{t=0}=p(\alpha).
$$

Let us change the time direction as $t\to-t$, then
\begin{equation}\label{t.44}
-\dot X=\nabla_\xi P(x,\Xi),
\end{equation}
$$
-\dot \Xi=-\nabla_x P(x,\Xi),.
$$
We want to solve the inverse problem:
\vspace {4mm}
$$
X(\alpha,0)=x(\alpha,T),\qquad \Xi(\alpha,0)=p(\alpha,T)
$$
The right-hand sides are considered as given data, and we are looking for
$X(\alpha,t)$, $\Xi(\alpha,t)$ for $0\le t\le T$. Obviously,
in our case the solutions have the form
$$
X(\alpha,t)=x(\alpha,T-t),\qquad \Xi(\alpha,t)=p(\alpha,T-t).
$$
Conclusion: we can use the "same" trajectories to move forward and backward in time.
But the incoming trajectories become outcoming and vice versa.

\begin{corollary}
Stable jumps become unstable.
\end{corollary}

But if there are no jumps (singularities of the projection mapping
$\pi:\Lambda_n^t\to\mathbb{R}^n_{x}$),
then our geometry (and the asymptotic solution!) is invertible in time.
This means that if we take the Cauchy problem solution $u$ for parabolic PDE such that
$$
u|_{t=0}=\exp{(-S_0(x)/\varepsilon)}\varphi^0(x),
$$
then the asymptotic solution for $t=T$ has the "WKB" form
$$
u_{\text{as}}|_{t=T}=\exp{(-S(x,T)/\varepsilon)}\varphi_0(x,T).
$$
Then taking the last function as the initial data for parabolic PDE in inverse time
(let $v_{\text{as}}(x,t)$ be its asymptotic solution), we get:
$$
v_{\text{as}}(x,T)=\exp{(-S_0(x)/\varepsilon)}\varphi^0(x)(1+0(\varepsilon))
$$
It can be easily verified in the case $P(x,\xi)=\xi^2$ i.e. for
simplest heat equation. If one constructs the solution of
inverse heat equation with initial data at $t=T$ of the form
$$
u_{\text{as}}|_{t=T}=\exp{(-S(x,T)/\varepsilon)}\varphi_0(x,T)
$$
using the Green function and calculates the integral be saddle
point method at $t=0$then the following result will be obtained
$$
\bigg(u_{\text{as}}|_{t=T}*G\bigg)|_{t=0}=\exp{(-S_0(x)/\varepsilon)}\varphi^0(x)(1+O(\varepsilon))
$$
We want to stress once again that this statement is true if
there is no singularities of projection mapping $\pi:
\Lambda_n^t\to \mathbb{R}^n_x$.

 But a jump brings problems:

\begin{figure}[h!] 
\includegraphics[width=0.7\textwidth]{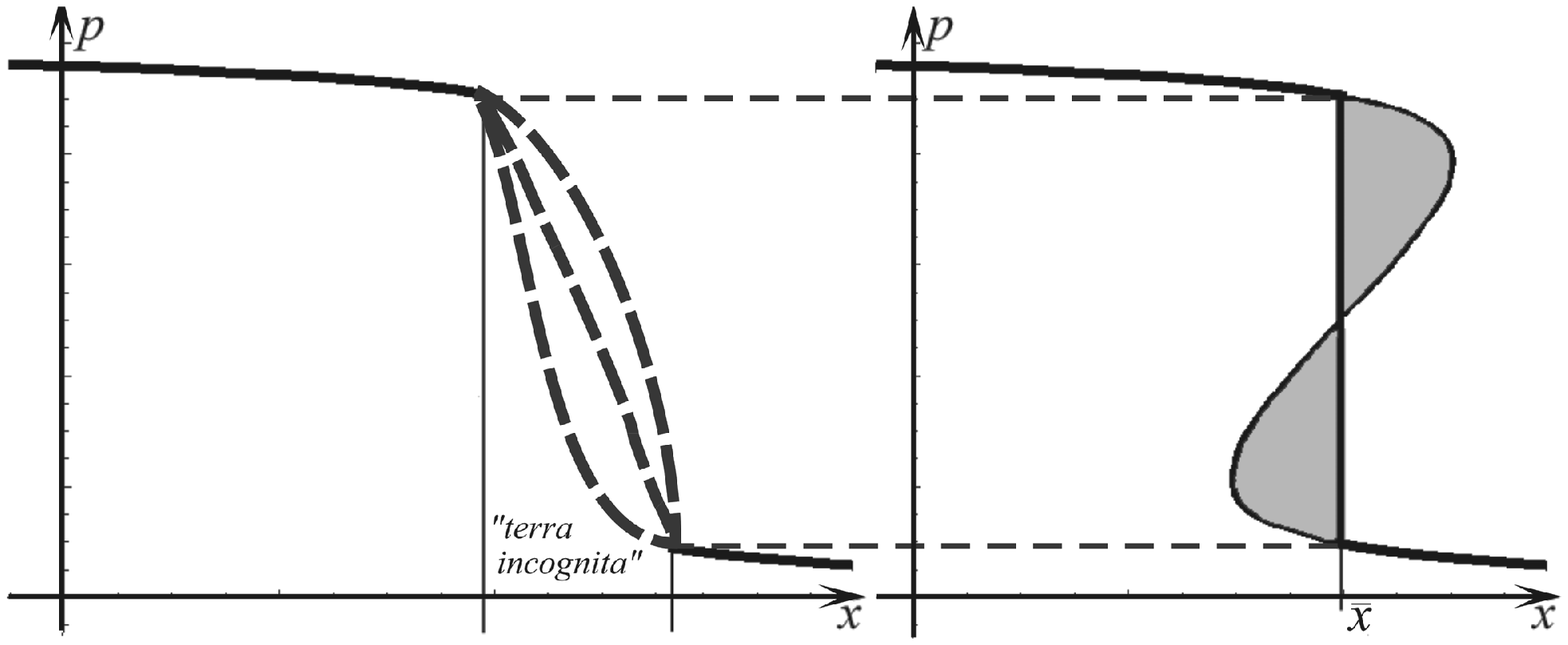}
\caption{}
\end{figure}

There is no unique reconstruction of the part of Lagrangian manifold
coming to the vertical line as $t$ increases (see Fig.9)!
But, fortunately, we can move ahead using the sense considerations.
The main point is that the function $S$ cannot attain its minimum (maximum)
inside the "terra incognita".
This allows one to calculate the integrals containing the reconstructed solution
without taking "terra incognita" into account
in the case where the integrand support contains this "terra incognita",
and we can formulate the following statement.

\begin{theorem}
Let the symbol $P(x,\xi)$ defined by ~\eqref{t.0} be such that the function
$B$ and measure $\mu$ do not depend on $x$ and $V=0$.
Let $u_{\varepsilon}(x,t)$ be a  solution of the Cauchy problem to
a Kolmogorov--Feller-type equation, and assume that,
for some $t\in (0,T)$,  there exists a logarithmic limit
$S(x,t)=\lim_{\varepsilon\to0}\big(-\varepsilon\ln u_{\varepsilon}\big)$,
namely,  the action function and the generalized amplitude
$\rho=w-\lim_{\varepsilon\to 0}\exp (2S/\varepsilon)u^2_{\varepsilon}$.

Let the singular support of $S(x,t)$ be a stratified manifold and $\rho_{\text{reg}}$
be the regular part of the generalized solution to the continuity equation
in the distribution sense.

Then, for an arbitrary test function $\phi=\phi(x)>0$
from the Schwartz space and for all $t\in[0,T]$,
the limit as ${\varepsilon}\to0$ of the integral
\begin{equation}\label{alpha1}
C\int\big(u_{\varepsilon}-\rho^{1/2}_{\text{reg}}
\exp(-S/{\varepsilon})\big)\phi dx
\end{equation}
is equal to~$0$, where $C^{-1}=\int u_{\varepsilon}\phi dx$.
\end{theorem}

This result can be extended to $C_0^{\infty}$-test functions and to the set of
smooth functions integrable with weight $\rho_{\text{reg}}^{1/2}\exp(-S/{\varepsilon})$
under the same assumption as in the theorem above.

Let $\Omega(t)$ be a subset in $R^n_{x}$, where we cannot define the functions
$S$ and $\varphi$
uniquely by the characteristics. Here we have two possible statements:
\begin{enumerate}

\item[(i)]
let $\Omega(t) \in \operatorname{supp}\phi$,
then relation (\ref{alpha1}) is true;

\item[(ii)]
assume that $x\in R^1$, $\Omega(t)$ is a union of segments $I_k$,
and the intersection between some $I_{\bar k}$ and
$\operatorname{supp}\phi$ is not empty but $I_{\bar k}$ does not belong to $\Omega(t)$,
then the limit of the integral as ${\varepsilon}\to0$
$$
\exp(\Psi/{\varepsilon})\int\big(u_{\varepsilon}-\rho^{1/2}_{\text{reg}}
\exp(-S/{\varepsilon})\big)\phi dx
$$
equals~$0$, where $\Psi$ is the minimal value of $S$ at the ends of $I_{\bar k}$.
\end{enumerate}

These statements actually mean that, from the viewpoint of the weak sense (momentum),
the density reconstructed arbitrarily inside the "terra incognita"
and according to the characteristics outside it
can be used in the same manner as the leading term of the asymptotic solution
constructed earlier by Maslov's tunnel canonical operator
and its modifications, \cite{M1,8,9,10,D5,DD}.

Now I will briefly speak about the proofs of the statements
about the invertibility in time that was formulated above.
Each of them can be divided into two parts.
First, it is to prove that, for all smooth reconstructions of Lagrangian manifold in a
"terra incognita" domain, the corresponding function $S_0=S_0(x)$ (for some fixed $t$)
cannot attain its minimum value inside this domain, see the lemma below.
This allows one to apply the Laplace method for calculating the integrals
mentioned in those statements taking into account that, due to this method and Lemma,
the results of these calculations do not depend on the values of the integrands
inside the "terra incognita" domain.
The proof is finished by taking account of the estimation
$$
u_{as}-u=O(\varepsilon)u_{as},
$$
which is true outside the singular support of the function $S$,
see Theorem above.
Now I formulate the Lemma.

\begin{lemma}
Let the symbol $P(x,\xi)$ defined by~\eqref{t.0} be such that the function $B$ and
the measure $\mu$ do not depend on $x$ and $V=0$.
Then the function $S$, i.e., the action function corresponding to the Lagrangian manifold,
cannot attain the minimal value inside the "terra incognita" domains.
\end{lemma}

We begin the consideration with a particular case when operator  symbol $P(x,\xi)$
does not depend on $x$
and restrict ourselves by one dimensional case studying. Let (a,b) is an interval
inside the "terra incognita" domain, and let $\bar x_0\in(a,b)$.
We will proceed by contradiction.
Assume that the function $S$ attain its minimal value at the point
$\bar x_0\in(a,b)$. We prove that, along the trajectory
of Hamilton system whose projection starts at $\bar x_0\in(a,b)$,
the following inequality is true:
\begin{equation}\label{d.7}
\frac{Dx}{Dx_0}|_{x_0=\bar x_0}\ne 0.
\end{equation}
This inequality leads to a contradiction because of the assumption that $\bar x_0$
belongs to the "terra incognita" domain, and hence it belongs to the projection
of the image of the singular (vertical) part of the Lagrangian manifold under backward in time shift
along the trajectories of the Hamilton system. In turn, this means that the projections of
all trajectories whose starting points are projected to the "terra incognita"
must intersect at a point for the forward in time motion.
So the above-mentioned inequality leads to a contradiction.
To prove this inequality, we write the projection of the Hamilton system trajectory
starting at $\bar x_0$. It has the form
\begin{equation}\label{d.15}
x=\bar x_0+tP_{\xi}(p_0),\qquad \qquad p_0=S_{0x_0}(\bar x_0).
\end{equation}
It is clear that $ p_0=S_{0x_0}(\bar x_0)=0$. Thus,
\begin{equation}\label{d.16}
\frac{\partial x}{\partial x_0}|_{x_0=\bar x_0}=1+tP_{\xi
\xi}(0)S_{0x_0x_0}(\bar x_0).
\end{equation}
Taking into account that that
\begin{equation}\label{d.8}
\frac{\partial^2 P}{\partial \xi^2}>0
\end{equation}
because of the convexity of $P$ and
\begin{equation}\label{d.99}
\frac{\partial^2 S_0}{\partial x^2_0}\ge 0
\end{equation}
due to the assumption that $\bar x_0$ is a point of minimal value,
we get the needed inequality ~\eqref{d.7}.
A multidimensional case differs from the case considered by changing the scalar values
$P_\xi$ and $S_{0x_0}$ by vectors (gradients).
In turn, this gives matrix inequalities in~\eqref{d.8} and~\eqref{d.99}.
The problem is to prove that the eigenvalues of the matrix
$P_{\xi \xi}(0)S_{0x_0x_0}(\bar x_0)$ are nonnegative
by using~\eqref{d.8} and~\eqref{d.99}.
For this, we can make a change of variables reducing the matrix $S_{0x_0x_0}(\bar x_0)$
to diagonal form.
This transformation induces the corresponding transformation
in the $p$-plane that transforms the matrix $P_{pp}(0)$
to a new symmetric positive matrix.
Now we note that the principal minors of the new matrix product
are products of matrix-factor principal minors
(because the second one is of diagonal form).
The determinants of the matrix-factor principal minors are nonnegative,
so the spectrum of the matrix $P_{\xi \xi}(0)S_{0x_0x_0}(\bar x_0)$
is also nonnegative and we again get~\eqref{d.7}.
To finish our consideration, we have to investigate the case where
the symbol $P=P(x,\xi)$ depends on $x$.
We again stay at the point $\bar x_0$, where the function $S_0$ attains its minimum.
If so, then
\begin{equation}\label{d.101}
p|_{t=0}= \frac{\partial S_0(x_0)}{\partial x_0}|_{x_0=\bar x_0}=0,
\end{equation}
and, by assumptions, $p=0$ for $t>0$.

The system for the matrices $\frac{\partial x}{\partial x_0}$ and
$ \frac{\partial p}{\partial x_0}$ follows from the Hamilton system
and has the form
\begin{equation}\label{d.111}
\frac{d}{dt}\frac{\partial x}{\partial x_0} =
\frac{\partial^2 P}{\partial
\xi \partial x}\frac{\partial x}{\partial x_0}+
\frac{\partial^2 P}{\partial \xi^2}\frac{\partial p}{\partial x_0}
\end{equation}
\begin{equation}\label{d.12}
\frac{d}{dt}\frac{\partial p}{\partial x_0} =
\frac{\partial^2 P}{\partial
x \partial \xi}\frac{\partial p}{\partial x_0}+
\frac{\partial^2 P}{\partial x^2}\frac{\partial x}{\partial x_0}
\end{equation}
Because of our assumptions (see the Lemma formulation),
we have
$P|_{\xi=0}=0$, $P_{xx}|_{\xi=0}=0$ and  $P_{\xi x}|_{\xi=0}=0$.
This means that, along the Hamilton system trajectory starting from
the point $x=\bar x_0$, $p=0$, we have
$p=0$ and $\frac{d}{dt}\frac{\partial p}{\partial x_0}=0$.
Thus, along the above-mentioned trajectory,
equations~\eqref{d.111} and~\eqref{d.12} have the form
\begin{equation}
\frac{d}{dt}\frac{\partial x}{\partial x_0} =
\frac{\partial^2 P(\bar x_0,0)}{\partial \xi^2}\frac{\partial p}{\partial x_0},
\end{equation}
\begin{equation} \label{d.21}
\frac{d}{dt}\frac{\partial p}{\partial x_0} =0.
\end{equation}
Integrating over $t$, we can transform the first equation to the form of Eq.~\eqref{d.16}
and apply all above arguments concerning this equality.

This note ends the proof.

The statement of the Lemma can be generalized as follows.
Let $B$ be a linear function in $x$,
$$
<B,\xi>=\sum_{k=1}^n b_k x_k\xi_k.
$$

Then, instead of \eqref{d.21}, we get
$$
\frac{d}{dt}\frac{\partial p}{\partial x_0}=\frac{\partial B}{\partial x}
\frac{\partial p}{\partial x_0}
$$
where the matrix $\frac{\partial B}{\partial x}$ has diagonal form
with elements equal to constants.

Then
$$
\frac{\partial p}{\partial x_0}=\exp^{t\frac{\partial B}{\partial x}}
\frac{\partial p_0}{\partial x_0}
$$
and
\begin{align*}
\frac{\partial x}{\partial x_0}
&=E+\frac{\partial^2 P}{\partial \xi^2}\int^t_0
\exp^{t'\frac{\partial p_0}{\partial x_0}}dt' \cdot
\frac{\partial p_0}{\partial x_0}
\\
&= \int_0^t\exp^{t'\frac{\partial B}{\partial x}}dt'
\bigg( \bigg(\int_0^t\exp^{t'\frac{\partial B}{\partial x}} dt'\bigg)^{-1}
+\frac{\partial^2 P}{\partial \xi^2}\frac{\partial p_0}{\partial x_0}\bigg).
\end{align*}
Here we used the fact that the matrix
$$
\int_0^t\exp^{t'\frac{\partial B}{\partial x}}dt'
$$
is in diagonal form with positive eigenvalues.

Thus we came to the relation with the same properties as \eqref{d.16}.
One can proceed further and generalize the statement to the case of an arbitrary drift.
But up to now the presence of the potential destroys our picture
and I will think about it.

\end{document}